\documentclass[a4paper,12pt]{article}
\usepackage{amssymb, amsmath}
\usepackage{graphics}
\usepackage[dvips]{graphicx}
\usepackage{eepic}
\usepackage{epic}
\usepackage{amsfonts}
\usepackage{latexsym}
\newtheorem{thm}{Theorem}[section]

\newtheorem{lmm}[thm]{Lemma}

\newtheorem{prp}[thm]{Proposition}

\newtheorem{question}[thm]{Question}
\newtheorem{observation}[thm]{Observation}

\newtheorem{dfn}[thm]{Definition}
\title{Wave Front Evolution and Pedal Evolution     
}
\author{Takashi \textsc{Nishimura}
\footnote{Research Group of Mathematical Sciences,  
Research Institute of Environment and Information Sciences, 
Yokohama National University, 
Yokohama 240-8501, JAPAN. 
\newline e-mail: \texttt{nishimura-takashi-yx@ynu.jp}}
}
\begin{document}
%

\maketitle

\begin{abstract}      
The calculus correspondence has been known to exist 
between generic pedal evolutions and generic wave front evolutions.      
In this paper, we first extend the known results on the calculus correspondence to evolutions with  multi-parameters,  
and then give applications of calculus correspondence.    Moreover, we discuss the 
possibility of generalization of the calculus correspondence to degenerate pedal evolutions and 
degenerate wave front evolutions.      
\end{abstract}
\section{Introduction}
\label{intro}
\noindent 
Throughout this paper, all maps, map-germs and vector fields are of class $C^\infty$ unless otherwise stated.   
\par 
A map-germ $\Phi: (\mathbb{R}^m,0)\to (\mathbb{R}^{m+1},0)$ is called 
a {{\it Legendrian map-germ}} if there exists a germ of unit vector field 
$\nu_\Phi$ 
along $\Phi$ 
such that the following 2 conditions hold, where the dot in the center stands for the scalar product of 
two vectors.      
\begin{enumerate}
\item 
$\frac{\partial \Phi}{\partial x_1}(x_1, \ldots, x_m)\cdot \nu_\Phi(x_1, \ldots, x_m)=\cdots =$ 
$\frac{\partial \Phi}{\partial x_m}(x_1, \ldots, x_m)\cdot \nu_\Phi(x_1, \ldots, x_m)$ 
$=0$.  
\item The map-germ $L_\Phi: (\mathbb{R}^{m},0)\to T_1\mathbb{R}^{m+1}$ defined by 
\[
L_\Phi(x_1, \ldots, x_m)=(\Phi(x_1, \ldots, x_m), \nu_\Phi(x_1, \ldots, x_m))
\] 
is non-singular, where $T_1\mathbb{R}^{m+1}$ is the unit tangent bundle of $\mathbb{R}^{m+1}$.      
\end{enumerate} 
The vector field $\nu_\Phi$, the map-germ $L_\Phi$ and the image of a Legendrian map-germ 
are called a {\it unit normal vector field of\/} $\Phi$, a {{\it Legendrian lift of\/} $\Phi$} and 
a {\it wave front\/} respectively.          
Singularities of Legendrian map-germs have been relatively well-studied 
(for instance, see \cite{arnoldetall, izumiya, zakalyukin2, zakalyukin}).  
\par 
\medskip 
Let ${\bf r}: (a,b)\to \mathbb{R}^2$ ($0\in (a,b)$) be a {non-singular plane curve without 
inflection point} (namely, a non-degenerate curve), and 
let $P$ be a point of $\mathbb{R}^2$.     
Then, the {{\it pedal curve of ${\bf r}$ 
{relative to the pedal point} $P$}} is defined as the trajectory of the foot of perpendicular to 
the tangent line $\{{\bf r}(s)+u{\bf r}'(s)\; |\; u\in \mathbb{R}\}$ at ${\bf r}(s)$ 
from $P$, and it is denoted by {$ped_{{\bf r}, P}: (a,b)\to \mathbb{R}^2$}.       The given point $P$ is called 
the {\it pedal point}.       
Let {$WF_{{\bf r}, P}: (a,b)\to \mathbb{R}^2$} be the solution curve of  
$$
\frac{d}{ds}WF_{{\bf r}, P}(s)=ped_{{\bf r}, P}(s)-P, \quad WF_{{\bf r}, P}(0)=(0,0), 
$$ 
where $s$ is the arc-length parameter of ${\bf r}$.     
Then, by definition of pedal curve, {${\bf r}'(s)$} (which is the unit tangent vector to ${\bf r}$ at 
${\bf r}(s)$) is a {unit normal vector to $WF_{{\bf r}, P}$ at $WF_{{\bf r}, P}(s)$}.    
Thus, the {Legendrian lift}  
$L_{WF_{{\bf r}, P}}: (a,b)\to T_1\mathbb{R}^2 $ 
given by 
$$
L_{WF_{{\bf r}, P}}(s)=\left(WF_{{\bf r}, P}(s), {\bf r}'(s)\right)
$$ 
is well-defined.      
Since the original curve ${\bf r}$ is without inflection point, 
by the Serret-Frenet formula (for the Serret-Frenet formula, see for instance \cite{brucegiblin}),  
$L_{WF_{{\bf r}, P}}$ is non-singular.        
Thus, the image of $WF_{{\bf r}, P}$ must be a {wave front curve}.        
\par 
Next, we move the pedal point $P$.       
Let $P: U\to \mathbb{R}^2$ be a map, 
where $U$ is an open neighborhood of the origin of $\mathbb{R}^n$,    
Then, we obtain two corank one maps 
$Un\mbox{-}ped_{{\bf r}, P}: (a,b)\times U\to (\mathbb{R}^2\times\mathbb{R}^{n},0)$ and 
$Un\mbox{-}WF_{{\bf r}, P}: (a,b)\times U\to (\mathbb{R}^2\times\mathbb{R}^{n},0)$ 
defined by 
$$
Un\mbox{-}ped_{{\bf r}, P}(s,u)=\left(ped_{{\bf r}, P(u)}(s), u\right)
\mbox{ and } 
Un\mbox{-}WF_{{\bf r}, P}(s,u)=\left(WF_{{\bf r}, P(u)}(s), u\right)  
$$
respectively.      
The map $Un\mbox{-}ped_{{\bf r}, P}$ (resp., $Un\mbox{-}WF_{{\bf r}, P}$) is called 
the {\it pedal unfolding of} $ped_{{\bf r}, P(0)}$ 
(resp., the {\it wave front unfolding of\/} $WF_{{\bf r}, P(0)}$).     
\par 
\medskip 
In \cite{arnold}, the evolution of generic wave fronts by time has been studied.         
We want to construct a different method from \cite{arnold} 
to study generic wave front evolutions.      
In order to do so, 
we pay attention to the relation between $Un\mbox{-}ped_{{\bf r}, P}$ and $Un\mbox{-}WF_{{\bf r}, P}$ since 
we have the following Proposition \ref{proposition 1} for $Un\mbox{-}ped_{{\bf r}, P}$.   
\begin{prp} \label{proposition 1}
Let ${\bf r}: (a,b)\to \mathbb{R}^2$ $(0\in (a,b))$ be a {non-singular plane curve without 
inflection point} such that ${\bf r}(0)=0$ and let $U$ be an  open neighborhood of the origin of $\mathbb{R}^n$. 
Moreover, we let $P: U\to \mathbb{R}^2$ be a map such that ${\bf r}(0)=P(0)=0$.    
Then, the map-germ  
$Un\mbox{-}ped_{{\bf r}, P}: ((a,b)\times U,(0,0))\to (\mathbb{R}^2\times\mathbb{R}^{n},(0,0))$ 
is $\mathcal{A}$-equivalent to  the normal form of {\rm (}Whitney umbrella{\rm )}$\times\mathbb{R}^{n-1}$ 
if and only if 
the origin $(0,0)$ of $(a,b)\times U$ is a regular 
point of the map $({\bf r},P): (a,b)\times U\to \mathbb{R}^{n+2}$ defined by 
$({\bf r},P)(s,u)=({\bf r}(s), P(u))$.   
\end{prp}
Here, two map-germs $f, g: (\mathbb{R}^{m},0)\to (\mathbb{R}^{m+1},0)$ are said to be 
$\mathcal{A}$-{\it equivalent\/} 
if there exist germs of diffeomorphism $h_s: (\mathbb{R}^{m},0) \to (\mathbb{R}^{m},0)$ 
and $h_t: (\mathbb{R}^{m+1},0) \to (\mathbb{R}^{m+1},0)$ 
such that $f=h_t\circ g\circ h_s$, and the 
{\it normal form of} ({\it Whitney umbrella})$\times\mathbb{R}^{n-1}$ is the map-germ defined by $(s,u)\mapsto (su_1, s^2,u)$ where 
$u=(u_1, \ldots, u_n)$.         
Proposition \ref{proposition 1} in the case $n=1$ is a special case of Theorem 1 in \cite{nishimura3}, 
Proof of Proposition \ref{proposition 1} is given in \S 2.     
\begin{figure}[hbtp]
\begin{center}
{\includegraphics[width=.3\linewidth]{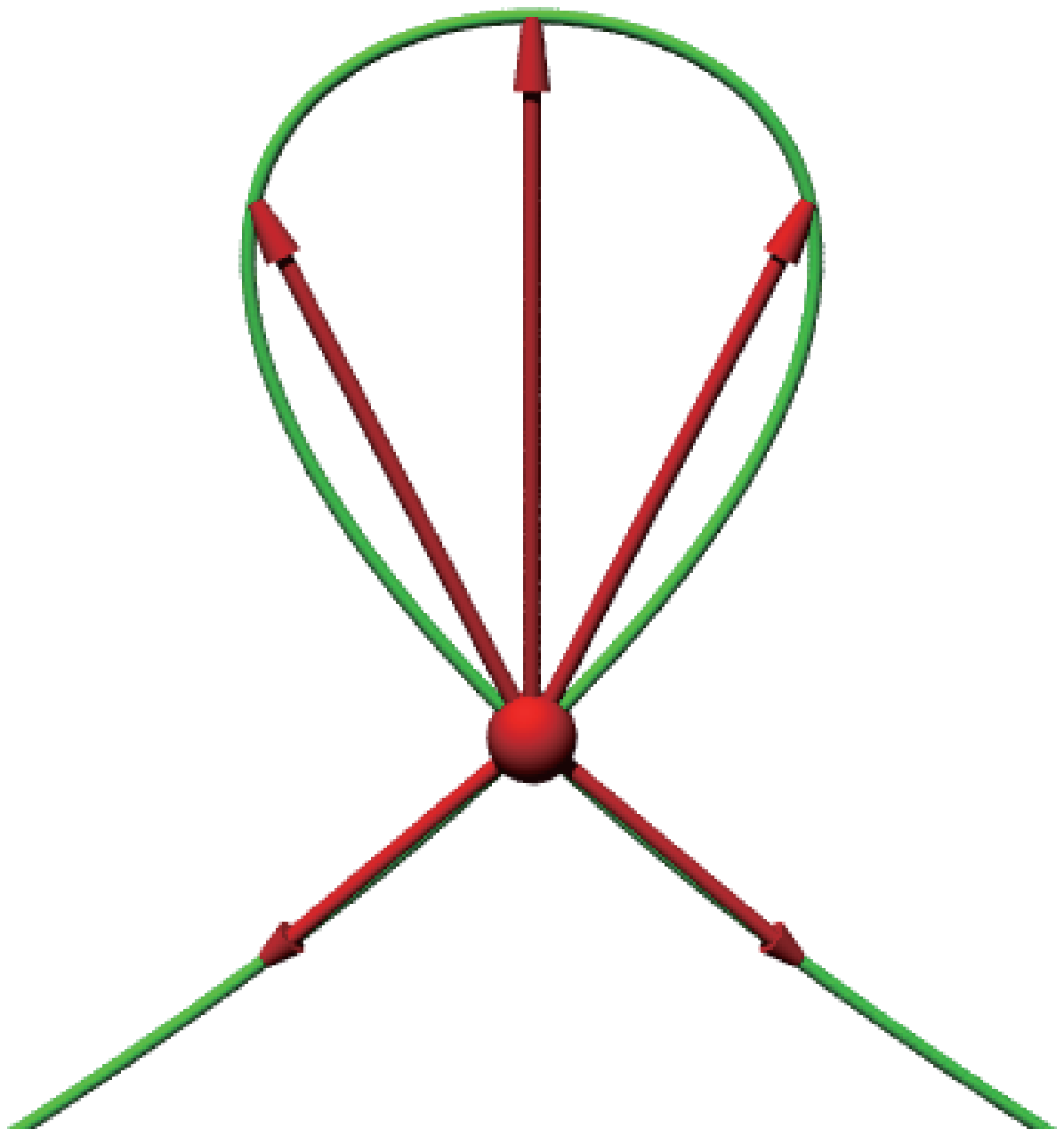}}
{\includegraphics[width=.3\linewidth]{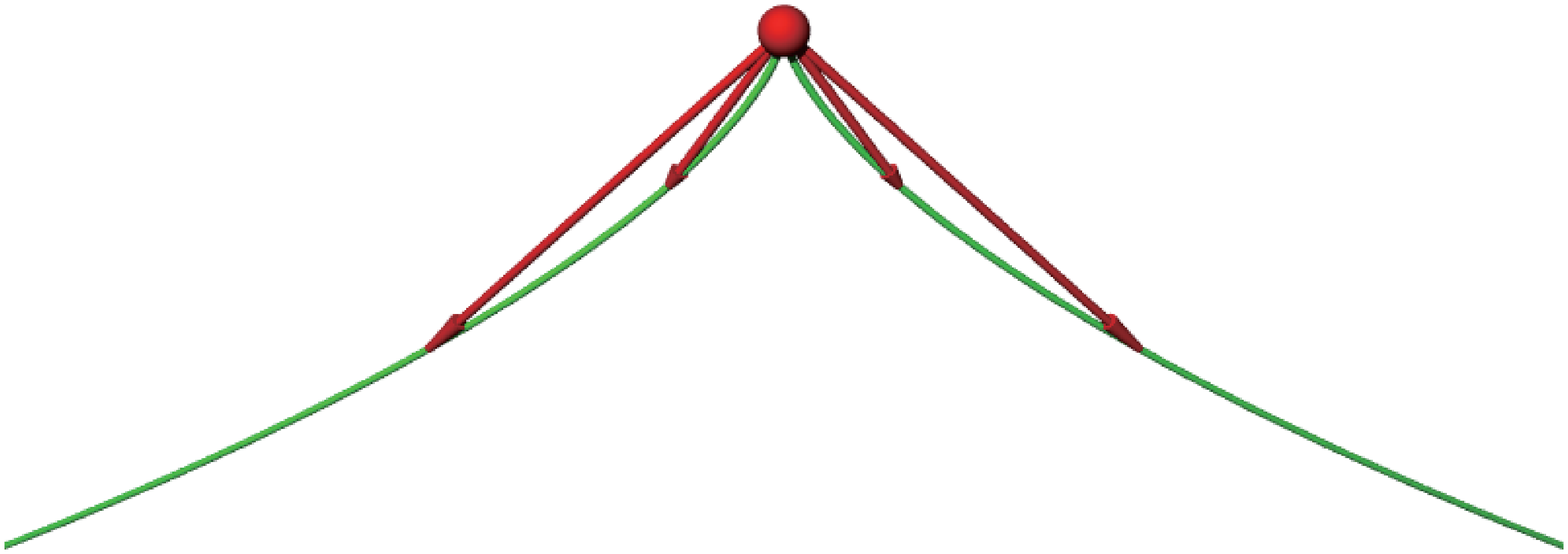}}
{\includegraphics[width=.3\linewidth]{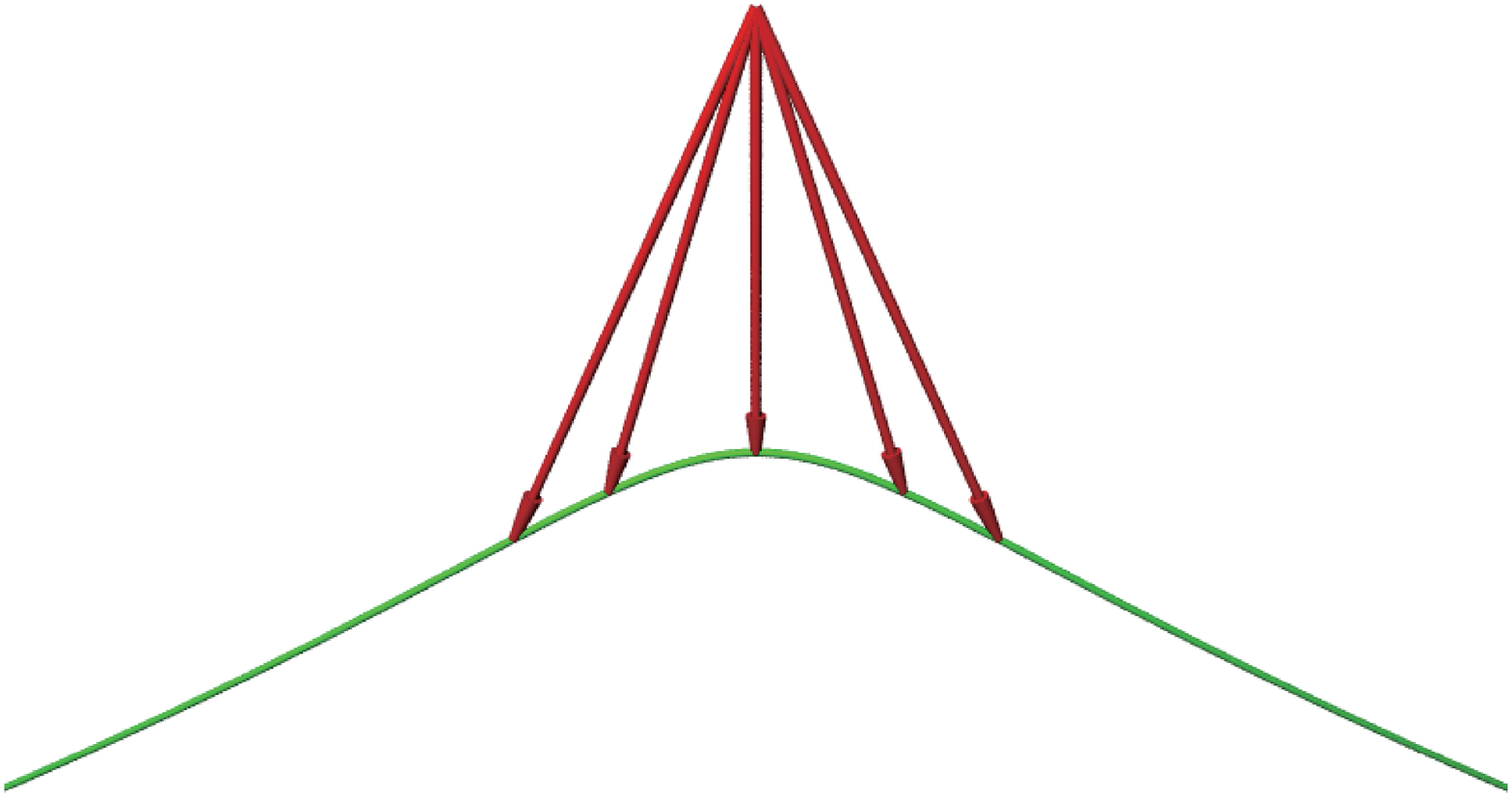}}
\caption{Pedal Evolution.}   
\label{figure 1}
\end{center}
\end{figure}   
\begin{figure}[hbtp]
\begin{center}
{\includegraphics[width=.3\linewidth]{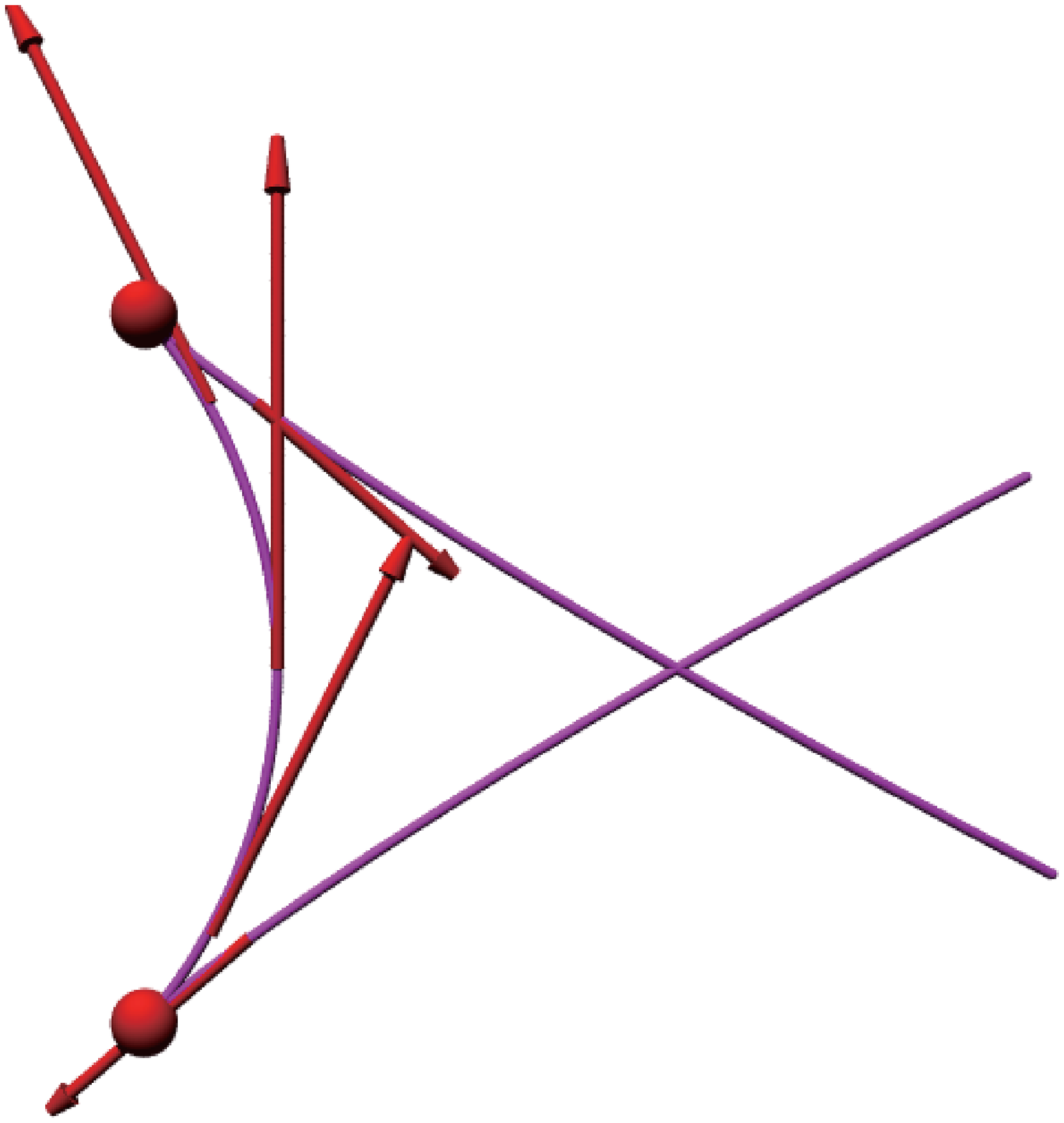}}
{\includegraphics[width=.3\linewidth]{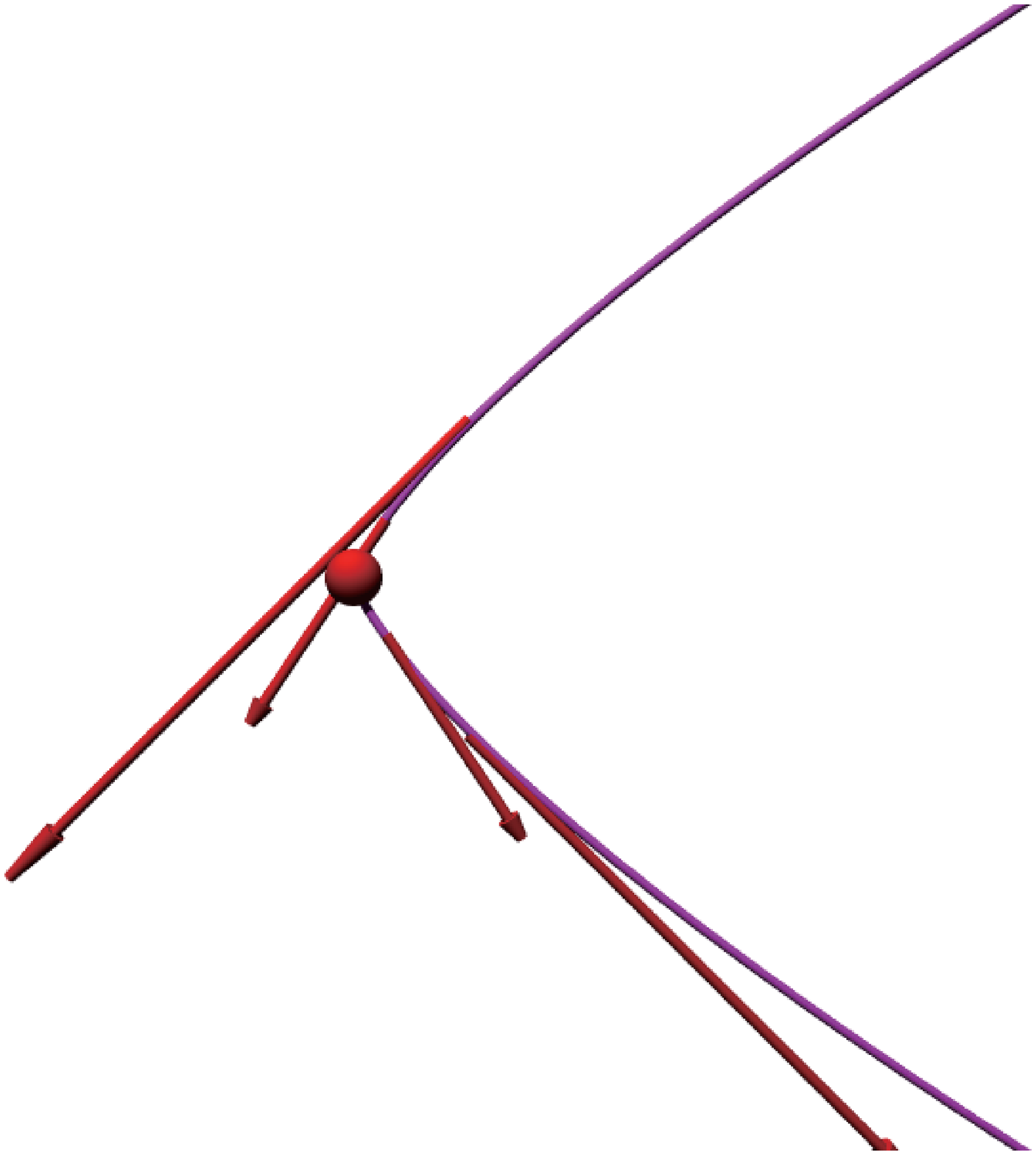}}
{\includegraphics[width=.3\linewidth]{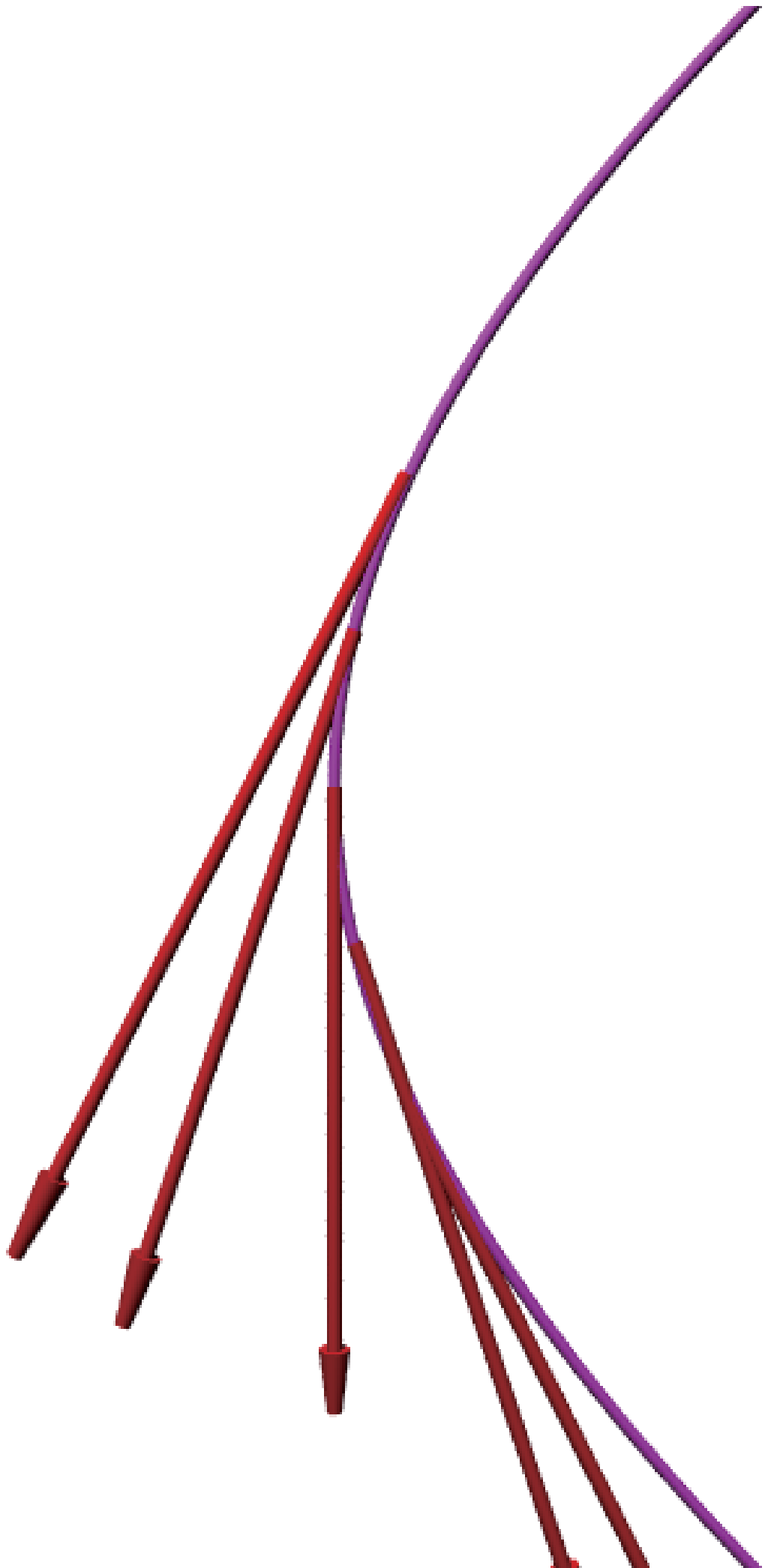}}
\caption{Wave Front Evolution.}  
\label{figure 2}
\end{center}
\end{figure}
\par 
By Figures \ref{figure 1} and \ref{figure 2}, it is easily conjectured that 
the pedal evolution $Un\mbox{-}ped_{{\bf r}, P}$ is $\mathcal{A}$-equivalent to 
the normal form of (Whitney umbrella)$\times\mathbb{R}^{n-1}$ 
if and only if the wave front evolution 
$Un\mbox{-}WF_{{\bf r}, P}$ is a (swallowtail)$\times \mathbb{R}^{n-1}$, 
where a { ({\it swallowtail\/})$\times\mathbb{R}^{n-1}$} is a map-germ $\mathcal{A}$-equivalent to 
$(s,u)\mapsto (3s^4+s^2u_1, -4s^3-2su_1, u)$ $(u=(u_1, \ldots, u_n))$; and  
in the case $n=1$ this conjecture has been actually proved in \cite{nishimurapacific} 
(such a correspondence is called the {\it calculus correspondence}.     
For more details on the known calculus correspondence, see \S 2).   
\par 
In this paper, 
we first extend the known results on the calculus correspondence to evolutions with multi-parameters,  
and then give applications of calculus correspondence.    Moreover, we discuss the 
possibility of generalization of calculus correspondence to degenerate pedal evolutions and 
degenerate wave front evolutions.    
%
%
\par 
\medskip 
In Section \ref{section 2}, known calculus correspondences are extended to evolutions with multi-parameters.   
The proof of Proposition \ref{proposition 1} is also given in Section \ref{section 2}.    
In Section \ref{section 3}, applications of calculus correspondence are given.    
Finally, the possibility of generalization of calculus correspondence is discussed in Section \ref{section 5}.    
\section{Extension of known calculus correspondences to evolutions with multi-parameters} \label{section 2}  
\begin{dfn} \label{definition 1}
{\rm 
A map-germ $\varphi: (\mathbb{R}\times \mathbb{R}^{n},(0,0))\to (\mathbb{R}^2\times\mathbb{R}^{n},(0,0))$ 
having the following form  
is said to be {\it of pedal unfolding type}.     
$$ \label{pedal}
\varphi(x,y)=\left(n(x,y)p(x,y), p(x,y), y\right)
$$
where $n: (\mathbb{R}\times\mathbb{R}^{n},(0,0))\to (\mathbb{R},0)$ is a function-germ 
satisfying $\frac{\partial n}{\partial x}(0,0)\ne 0$, 
$p: (\mathbb{R}\times\mathbb{R}^{n},(0,0))\to (\mathbb{R},0)$ is a function-germ and $y=(y_1, \ldots, y_n)$.    
}
\end{dfn}
In the case $n=1$, Definition \ref{definition 1} has been given in \cite{nishimurapacific}. 
\begin{prp}\label{proposition 2} 
Let ${\bf r}: (a,b)\to \mathbb{R}^2$ $(0\in (a,b))$ be a non-singular plane curve without 
inflection point such that ${\bf r}(0)=0$ and 
let $P: U\to \mathbb{R}^2$ be a map such that $P(0)={\bf r}(0)=0$, where $U$ is an open neighborhood of the origin of $\mathbb{R}^n$.     
Then, $Un\mbox{-}ped_{{\bf r}, P}$ is $\mathcal{A}$-equivalent to a map-germ of pedal unfolding type.   
\end{prp}
\textit{Proof of Proposition \ref{proposition 2}.}\qquad 
Since ${\bf r}: (a,b)\to \mathbb{R}^2$ $(0\in (a,b))$ is  
non-singular, 
we may assume that ${\bf r}(x)=(-x, r_2(x))$, $r_2(0)=\frac{d r_2}{dx}(0)=0$ near $0$.    
Put {$n(x)=\frac{d r_2}{dx}(x)$}.      
Then, $(n(x), 1)$ is a normal vector to ${\bf r}$ at 
${\bf r}(x)$.     Since ${\bf r}$ is without inflection point, 
we have that {$\frac{d{n}}{dx}(x)\ne 0$}.    
Then, since ${\bf r}(0)=P(0)=0$, 
$ped_{{\bf r},P(0)}(x)\in T_{P(0)}\mathbb{R}^2$ has the form:
$$
{ped_{{\bf r},P(0)}(x)}=p(x)(n(x), 1){=\left(n(x)p(x), p(x)\right)}.    
$$
Therefore, we have:
\begin{eqnarray*}
Un\mbox{-}ped_{{\bf r}, P}(x,y) & = & \left({ped_{{\bf r},P(y)}(x)}, y\right) \\ 
{ } & = & 
\left(p(x,y)\left(n(x), 1\right)+P(y), y\right) \\ 
{ } & = & \left(\left(n(x)p(x,y), p(x,y)\right)+P(y), y\right) \\ 
{ } & \sim_{\mathcal{A}} & \left((n(x)p(x,y), p(x,y), y\right) .    
\end{eqnarray*}
\hfill\qquad $\Box$ 
\par 
\medskip 
\begin{prp}\label{proposition 3}
Let ${\bf r}: (a,b)\to \mathbb{R}^2$ $(0\in (a,b))$ be a non-singular plane curve without 
inflection point such that ${\bf r}(0)=0$ and 
let $P: U\to \mathbb{R}^2$ be a map such that $P(0)={\bf r}(0)=0$, where $U$ is an open neighborhood of the origin of $\mathbb{R}^n$.     
Then, $Un\mbox{-}ped_{{\bf r}, P}$ is $\mathcal{A}$-equivalent to a map-germ of the form:   
$(x,y)\mapsto \left(x(x^2+q(y)), x^2+q(y), y\right)$.   
\end{prp}
\textit{Proof of Proposition \ref{proposition 3}.}\qquad 
As same as the proof of Proposition \ref{proposition 2}, 
we may assume that ${\bf r}(x)=(-x, r_2(x))$, $r_2(0)=\frac{d r_2}{dx}(0)=0$ near $0$.       
By definition of pedal curve, the following holds:
$$
n(x)\left(n(x)p(x)+x\right)+\left(p(x)-r_2(x)\right)=0.
$$  
Here, $n(x)$ and $p(x)$ are functions defined in the proof of Proposition \ref{proposition 2}.     
Thus, we have the following locally:
$$
p(x)=\frac{r_2(x)-x n(x)}{1+n^2(x)}.
$$
Since ${\bf r}$ is without inflection point and $n(x)=\frac{d r_2}{dx}(x)$, 
there exists a function $\xi(x)$ such that 
$r_2(x)-x n(x)=x^2\xi(x)$ and $\xi(0)\ne 0$ by Hadamard's lemma 
(for Hadamard's lemma, see \cite{brucegiblin}).     
Thus, $p: (\mathbb{R},0)\to (\mathbb{R},0)$ is a Morse function-germ.     
By the Morse lemma with parameters (see \cite{brucegiblin}), we have that 
$$
{Un\mbox{-}ped_{{\bf r},P}(x,y)}  {=}  
{\left(ped_{{\bf r},P(y)}(x),y\right)} 
$$ 
is $\mathcal{A}$-equivalent to 
a map-germ of the form: 
\begin{equation}\label{equation 1}
(x,y)\mapsto {\left(n(x,y)\left(x^2+q(y)\right), x^2+q(y), y\right)}. 
\end{equation}
Since $\frac{\partial n}{\partial x}(0)\ne 0$, 
by using the Malgrange preparation theorem (for the Malgrange preparation theorem, 
for instance see \cite{arnoldetall}), 
$(x,y)\mapsto {\left(n(x,y)\left(x^2+q(y)\right), x^2+q(y), y\right)}$ is $\mathcal{A}$-equivalent to 
\begin{equation}\label{equation 2}
(x,y)\mapsto \left(x\left(x^2+q(y)\right), x^2+q(y), y\right).
\end{equation}
\hfill\qquad $\Box$ 
\par 
\noindent 
\smallskip 
Note that Proposition \ref{proposition 3} may be proved without using the Malgrange preparation theorem.  
Alternatively, we may adopt a simple method used to prove the criterion of cuspidal crosscap given in \cite{fujimorisajiumeharayamada}.        
Namely, by dividing $n(x,y)$ into the sum of an odd function and an even function with respect to     
the variable $x$, it is possible to show that the map-germ (\ref{equation 1})
is $\mathcal{A}$-equivalent to the map-germ (\ref{equation 2}). 
\begin{dfn}[\cite{mond}]\label{definition S_k}
{\rm Let $k$ be a non-negative integer.   
Then, a map-germ $f: (\mathbb{R}^2,0)\to (\mathbb{R}^3,0)$ is said to be {\it of $S_k$ type} 
if $f$ is $\mathcal{A}$-equivalent to the map-germ $(x,y)\mapsto (x(x^2\pm y^{k+1}), x^2\pm y^{k+1}, y)$.   
}
\end{dfn}
\noindent 
By Proposition \ref{proposition 3}, 
singularities of one-parameter pedal unfoldings $Un\mbox{-}ped_{{\bf r}, P}$ must be 
of $S_k$ type ($k$ is a certain non-zero integer) 
provided that ${\bf r}$ is non-degenerate, ${\bf r}(0)=P(0)$ and $q(y)$ is not flat.        
This fact has been already proved in \cite{nishimura3} by using a characterization of spherical pedal given in 
\cite{nishimura}.        Thus, the proof given here is an alternative proof.             
Even if we moved the given non-degenerate curve ${\bf r}$ depending on the parameter $y$, 
the proof of Proposition \ref{proposition 3} 
shows that new singularities never occur for the map-germ of the form 
$(x,y)\mapsto \left(ped_{{\bf r}_y,P(y)}(x), y\right)$ 
provided that ${\bf r}_0$ is non-degenerate, ${\bf r}_0(0)=P(0)$ and $q(y)$ is not flat.       
\par 
\bigskip 
\textit{Proof of Proposition \ref{proposition 1}.}\footnote{The author's original proof of Proposition 
\ref{proposition 1} used Mather's infinitesimal characterization of stable map-germs (\cite{mather2})) and 
Mather's classification theorem (\cite{mather4}).    
The proof given here, which is self-contained, was suggested by the referee.   
}
\qquad 
By Proposition \ref{proposition 3}, $Un\mbox{-}ped_{{\bf r},P}$ is $\mathcal{A}$-equivalent to 
a map-germ $\psi(x,y)=\left(x(x^2+q(y)), x^2+q(y), y\right)$ under the assumption of Proposition \ref{proposition 1}.         
It is easily seen that 
the origin $(0,0)$ of $(a,b)\times U$ is a regular 
point of the map $({\bf r},P): (a,b)\times U\to \mathbb{R}^{n+2}$
if and only if 
there exists an integer $i$ $(1\le i\le n)$ 
such that $\frac{\partial q}{\partial y_i}(0)\ne 0$ for $q(y)$.     
Thus, it is sufficient to show that 
$\psi(x,y)=\left(x(x^2+q(y)), x^2+q(y), y\right)$ 
is a (Whitney umbrella)$\times\mathbb{R}^{n-1}$ if and only if 
there exists an integer $i$ $(1\le i\le n)$ 
such that $\frac{\partial q}{\partial y_i}(0)\ne 0$ for $q(y)$.    
\par 
Suppose that $\psi$ is a (Whitney umbrella)$\times\mathbb{R}^{n-1}$.     
Let $S_1\subset J^2(\mathbb{R}^{n+1}, \mathbb{R}^{n+2})$ be the set of corank one $2$-jets.      
Note that $j^2f$ is transverse to $S_1$ where $f$ denotes 
the normal form of (Whitney umbrella)$\times\mathbb{R}^{n-1}$.      
Since $\psi$ is $\mathcal{A}$-equivalent to $f$, $j^2\psi$, too, is transverse to $S_1$.    
Since $S_1(\psi)=(j^2\psi)^{-1}(S_1)=\{3x^2+q(y)=0, x=0\}$, $j^2\psi$ is transverse to $S_1$ 
if and only if 
$$
\mbox{\rm rank}\left( d(3x^2+q(y)), dx\right)=\mbox{\rm rank}
\left(
\begin{array}{cc}
6 x & 1 \\
d q(y) & 0
\end{array}
\right)
= 2.    
$$  
Therefore, 
there exists an integer $i$ $(1\le i\le n)$ 
such that $\frac{\partial q}{\partial y_i}(0)\ne 0$. 
\par  
Conversely, suppose that 
there exists an integer $i$ $(1\le i\le n)$ 
such that $\frac{\partial q}{\partial y_i}(0)\ne 0$. 
Set 
\begin{eqnarray*}
h_s(x, y_1, \ldots, y_n) & = & (x, y_1, \ldots, y_{i-1}, x^2+q(y), y_{i+1}, \ldots, y_n), \\  
H_t(X_1, X_2, Y_1, \ldots, Y_n) & = & (X_1, Y_i, Y_1, \ldots, Y_{i-1}, -X_2+Y_i, Y_{i+1}, \ldots, Y_n)). 
\end{eqnarray*}     
Then, $h_s$ (resp., $H_t$) is a germ of diffeomorphism of $(\mathbb{R}\times\mathbb{R}^n, (0,0))$ 
(resp., $(\mathbb{R}^2\times\mathbb{R}^n, (0,0))$).         
Set also $f_i(x, y_1, \ldots, y_n)=(xy_i, x^2, y_1, \ldots, y_n)$.      
Then, we have: 
$$
H_t\circ f_i\circ h_s\left(x, y_{1}, \ldots, y_n\right)=
\left(x(x^2+q(y)), x^2+q(y), y_1, \ldots, y_{i-1}, q(y), y_{i+1}, \ldots, y_n\right).   
$$
Since $\frac{\partial q}{\partial y_i}(0)\ne 0$, the map-germ 
$(y_1, \ldots, y_n)\mapsto (y_1, \ldots, y_{i-1}, q(y), y_{i+1}, \ldots, y_n)$ is a germ of diffeomorphism.    
Thus, $H_t\circ f_i\circ h_s$ is $\mathcal{A}$-equivalent to $\psi$.     
Since the map-germ $f_i$ is clearly a (Whitney umbrella)$\times\mathbb{R}^{n-1}$, 
$\psi$ must be a (Whitney umbrella)$\times\mathbb{R}^{n-1}$.    
\hfill\qquad $\Box$ 
\begin{dfn} \label{definition 2}
{\rm 
For a map-germ of pedal unfolding type 
$$
\varphi(x,y)=(n(x,y)p(x,y), p(x,y),y), 
$$ 
set  
{
\[
\mathcal{I}(\varphi)(x,y)=\left(
\int_0^x n(x,y)p(x,y)dx, \int_0^x p(x,y)dx, y
\right). 
\]
}
The map-germ $\mathcal{I}(\varphi): (\mathbb{R}\times\mathbb{R}^n,0)\to (\mathbb{R}^2\times\mathbb{R}^n,0)$ is 
called the {{\it integration of }$\varphi$}.    
}
\end{dfn}
In the case $n=1$, Definition \ref{definition 2} has been given in \cite{nishimurapacific}. 
\begin{dfn}\label{normalized}
{\rm 
A Legendrian map-germ 
$\Phi:  (\mathbb{R}\times\mathbb{R}^n,(0,0))\to (\mathbb{R}^2\times\mathbb{R}^n,(0,0))$ 
is said to be {\it normalized\/} if $\Phi$ satisfies the following three conditions: 
\begin{enumerate}
\item The map-germ $\Phi$ has the form $\Phi(x,y)=(\Phi_1(x,y), \Phi_2(x,y),y)$ where $y=(y_1, \ldots, y_n)$.  
\item The condition $\frac{\partial \Phi_2}{\partial x}(0,0)=0$ holds.    
\item The vector $\nu_\Phi(0,0)$ is $\frac{\partial}{\partial X_1}$ or $-\frac{\partial}{\partial X_1}$, 
where $(X_1, X_2, Y_1, \ldots, Y_n)$ denotes the standard coordinate system of 
$(\mathbb{R}^2\times\mathbb{R}^n, (0,0))$.     
\end{enumerate}    
}
\end{dfn}
In the case $n=1$, Definition \ref{normalized} has been given in \cite{nishimurapacific}. 
\begin{dfn}\label{definition 5}
{\rm 
For a normalized Legendrian map-germ  
$$\Phi(x,y)=(\Phi_1(x,y), \Phi_2(x,y),y),$$ 
set  
{
\[
\mathcal{D}(\Phi)(x,y)=\left(
\frac{\partial \Phi_1}{\partial x}(x,y), 
\frac{\partial \Phi_2}{\partial x}(x,y), y
\right). 
\]
}
The map-germ $\mathcal{D}(\Phi): (\mathbb{R}\times\mathbb{R}^n,(0,0))\to 
(\mathbb{R}^2\times\mathbb{R}^n,(0,0))$ 
is called the {{\it differential of }$\Phi$}.        
}
\end{dfn}
In the case $n=1$, Definition \ref{definition 5} has been given in \cite{nishimurapacific}. 
\begin{prp}\label{proposition 4}
\begin{enumerate}
\item For a map-germ of pedal unfolding type 
$\varphi : (\mathbb{R}\times\mathbb{R}^n,(0,0))\to 
(\mathbb{R}^2\times\mathbb{R}^n,(0,0))$, 
$\mathcal{I}(\varphi)$ is a normalized Legendrian map-germ.  
\item For a normalized Legendrian map-germ  
$\Phi : (\mathbb{R}\times\mathbb{R}^n,(0,0))\to 
(\mathbb{R}^2\times\mathbb{R}^n,(0,0))$, 
$\mathcal{D}(\Phi)$ is a map-germ of pedal unfolding type.  
\end{enumerate}
\end{prp}
\noindent 
In the case $n=1$, 
Proposition \ref{proposition 4} with its proof can be found in \cite{nishimurapacific}.     
The proof given in \cite{nishimurapacific} works well even to the case $n\ge 2$.     
\par 
\medskip 
The following set is denoted by $\mathcal{W}$.    
$$ 
{\footnotesize \left\{ \varphi: (\mathbb{R}\times\mathbb{R}^n,(0,0))\to (\mathbb{R}^2\times\mathbb{R}^n,(0,0))
\; \mbox{ ({W}hitney umbrella)$\times\mathbb{R}^{n-1}$, pedal unfolding type}\right\}.} 
$$
And set also 
\begin{eqnarray*} 
{\mathcal{S}}  & = &   
\left\{ \Phi: (\mathbb{R}\times\mathbb{R}^n,(0,0))\to (\mathbb{R}^2\times\mathbb{R}^n,(0,0))\; 
\mbox{ normalized (swallowtail)$\times\mathbb{R}^{n-1}$}
\right\}, \\  
{\mathcal{N}}  
& = & 
{\small\small \left\{ \varphi: (\mathbb{R}\times\mathbb{R}^n,(0,0))\to (\mathbb{R}^2\times\mathbb{R}^n,(0,0))\; 
\mbox{ {non-singular, pedal unfolding type}}\right\},} \\ 
{\mathcal{C}}  & = &   
\left\{ \Phi: (\mathbb{R}\times\mathbb{R}^n,(0,0))\to (\mathbb{R}^2\times\mathbb{R}^n,(0,0))\; 
\mbox{ normalized (cusp)$\times\mathbb{R}^{n}$}
\right\}, 
\end{eqnarray*}
where a map-germ $\Phi: (\mathbb{R}\times\mathbb{R}^n,(0,0))\to (\mathbb{R}^2\times\mathbb{R}^n,(0,0))$ 
is called a ({\it cusp\/})$\times\mathbb{R}^n$ 
if it is $\mathcal{A}$-equivalent to $(x,y)\mapsto (2x^3,-3x^2,y)$ $(y=(y_1, \ldots, y_n))$.               
The following Theorems \ref{theorem 1} and \ref{theorem 2} 
are extensions of known calculus correspondences to multi-parameters.   
\begin{thm} \label{theorem 1}
\begin{enumerate}
\item The map $\mathcal{I}: \mathcal{W}\to \mathcal{S}$ 
defined by $\mathcal{W}\ni\varphi\mapsto \mathcal{I}(\varphi)\in\mathcal{S}$ 
is well-defined and bijective. 
\item The map $\mathcal{D}: \mathcal{S}\to \mathcal{W}$ 
defined by $\mathcal{S}\ni\Phi\mapsto \mathcal{D}(\Phi)\in\mathcal{W}$ 
is well-defined and bijective. 
\end{enumerate}     
\end{thm}
\begin{thm} \label{theorem 2}
\begin{enumerate}
\item The map $\mathcal{I}: \mathcal{N}\to \mathcal{C}$ 
defined by $\mathcal{N}\ni\varphi\mapsto \mathcal{I}(\varphi)\in\mathcal{C}$ 
is well-defined and bijective.    
\item The map $\mathcal{D}: \mathcal{C}\to \mathcal{N}$ 
defined by $\mathcal{C}\ni\Phi\mapsto \mathcal{D}(\Phi)\in\mathcal{N}$ 
is well-defined and bijective.    
\end{enumerate}     
\end{thm}
\noindent 
In the case $n=1$, the proofs of Theorems \ref{theorem 1} and \ref{theorem 2} 
can be found in \cite{nishimurapacific}.     
For the proof of Theorem \ref{theorem 1} in the case $n=1$,  two criteria 
(Theorems \ref{theorem 3} and \ref{theorem 4}) have been used in \cite{nishimurapacific}.     
Theorem \ref{theorem 4} works well even in the case $n\ge 2$.       
Although it is uncertain that Theorem \ref{theorem 3} works well 
even in the case $n\ge 2$, since a (Whitney umbrella)$\times\mathbb{R}^{n-1}$ is stable, 
by Mather's infinitesimal characterization of stable map-germs,  
Theorem \ref{theorem 1} in general case can be proved.         
On the other hand, for the proof of Theorem \ref{theorem 2} in the case $n=1$,  
Theorem \ref{theorem 3} has not been used in \cite{nishimurapacific} though 
Theorem \ref{theorem 4} has been used.       
Hence the proof of Theorem \ref{theorem 2} works well even in the case $n\ge 2$.  
\par 
Besides Theorems \ref{theorem 1} and \ref{theorem 2}, 
there is one more example of calculus correspondence (Proposition \ref{proposition 5}).    
Since Proposition \ref{proposition 5} is almost trivial, its proof is omitted.      
Put 
\begin{eqnarray*}
{ } & { } & {\small \mathcal{N}_{\mbox{non-zero}}} \\ 
{ } & = & 
{\small \left\{ \varphi: 
(\mathbb{R}\times\mathbb{R}^n,(0,0))\to \mathbb{R}^2\times\mathbb{R}^n\mbox{-}\{ (0, 0)\} \; 
\mbox{non-singular, of pedal unfolding type} \right\} , }\\ 
{ } & { } & \widetilde{\mathcal{N}} \\ 
{ } & = & 
\left\{ \Phi: (\mathbb{R}\times\mathbb{R}^n,(0,0))\to (\mathbb{R}^2\times\mathbb{R}^n,(0,0))\; 
\mbox{ normalized non-singular Legendrian}
\right\}.  
\end{eqnarray*}
\begin{prp} \label{proposition 5}
\begin{enumerate}
\item The map $\mathcal{I}: \mathcal{N}_{\mbox{non-zero}}\to \widetilde{\mathcal{N}}$ 
defined by $\mathcal{N}_{\mbox{non-zero}}\ni\varphi\mapsto \mathcal{I}(\varphi)\in\widetilde{\mathcal{N}}$ 
is well-defined and bijective.    
\item The map $\mathcal{D}: \widetilde{\mathcal{N}}\to \mathcal{N}_{\mbox{non-zero}}$ 
defined by $\widetilde{\mathcal{N}}\ni\Phi\mapsto \mathcal{D}(\Phi)\in\mathcal{N}_{\mbox{non-zero}}$ 
is well-defined and bijective.    
\end{enumerate}
\end{prp}   
\begin{dfn}[\cite{mond}]
{\rm Let $T: \mathbb{R}^2\to \mathbb{R}^2$ be the linear transformation of the form 
$T(s,\lambda)=(-s,\lambda)$.    
Two function germs $p_1, p_2: (\mathbb{R}^2,0)\to (\mathbb{R},0)$ are said to be 
$\mathcal{K}^T$-{\it equivalent\/} if there exists a germ of diffeomorphism 
$h : (\mathbb{R}^2,0)\to (\mathbb{R}^2,0)$ having the form $h\circ T=T\circ h$ and 
a function-germ   
$M : (\mathbb{R}^2,(0,0))\to \mathbb{R}$ having the form $M\circ T=M$, $M(0,0)\ne 0$ 
such that $p_1\circ h(s,\lambda)=M(s,\lambda)p_2(s,\lambda)$.    
}
\end{dfn}
\begin{thm}[\cite{mond}]   \label{theorem 3}
Two map-germs $f_i: (\mathbb{R}^2,0)\to (\mathbb{R}^3,0)$ $(i=1,2)$ of the following form 
$$
f_i(x, y)=(n_i(x,y)p_i(x^2, y),x^2,y),
$$
where 
$\frac{\partial n_i}{\partial x}(0,0)\ne 0$ and $p_i(x,y)$ is not flat for each $i\in \{1,2\}$,  
are $\mathcal{A}$-equivalent if and only if the function-germs $p_i(x^2, y)$ are 
$\mathcal{K}^T$-equivalent.   
\end{thm}
\begin{dfn}\label{definition 3}
{\rm 
Let $\Phi: (\mathbb{R}\times\mathbb{R}^n, (0,0))\to 
(\mathbb{R}^2\times\mathbb{R}^n,(0,0))$ be a Legendrian map-germ and let 
$\nu_\Phi$ be a unit normal vector field of $\Phi$ given in the definition of Legendrian map-germs.      
The function-germ $LJ_{\Phi}: (\mathbb{R}\times\mathbb{R}^n,0)\to \mathbb{R}$ defined by the following is 
called the {{\it Legendrian-Jacobian} of $\Phi$} where $(x,y)=(x, y_1, \ldots, y_n)$.    
\[
LJ_{\Phi}(x,y)= 
\det 
\left(
\frac{\partial \Phi}{\partial x}(x,y), \frac{\partial \Phi}{\partial y_1}(x,y), 
\ldots, \frac{\partial \Phi}{\partial y_n}(x,y), {\nu}_\Phi(x,y)
\right).    
\]
}
\end{dfn}
In the case $n=1$, Definition \ref{definition 3} can be found in \cite{nishimurapacific}.    
Note that if $\nu_\Phi$ satisfies the conditions of unit normal vector field of $\Phi$, 
then $-\nu_\Phi$ also satisfies 
them.   
Thus, the sign of $LJ_{\Phi}(x,y)$ depends on the particular choice of unit normal vector field $\nu_\Phi$.     
The Legendrian Jacobian of $\Phi$ is called also 
the {\it signed area density function} 
(for instance, see \cite{sajiumeharayamada2}).       
Although it seems reasonable to call $LJ_{\Phi}$ the area density function from the viewpoint of 
investigating the singular surface $\Phi(U)$ ($U$ 
is a sufficiently small neighborhood of the origin of $\mathbb{R}^2$), 
it seems reasonable to call it the Legendrian Jacobian from the viewpoint of 
investigating the singular map-germ $\Phi$.      
\begin{thm}[\cite{sajiumeharayamada}] \label{theorem 4}
Let  
$\Phi: (\mathbb{R}\times\mathbb{R}^n,(0,0))\to (\mathbb{R}^2\times\mathbb{R}^n,(0,0))$ 
be a normalized Legendrian map-germ, 
\begin{enumerate}
\item  $\Phi$ is a {\rm (}swallowtail\/{\rm )}$\times\mathbb{R}^{n-1}$  
if and only if  the following two hold where $y=(y_1, \ldots y_n)$:
$$
Q\left(LJ_\Phi, \frac{\partial LJ_\Phi}{\partial x}\right)\cong Q(x,y_1),\; 
\frac{\partial^2 LJ_\Phi}{\partial x^2}(0,0)\ne 0.     
$$   
\item 
$\Phi$ is a {\rm (}cusp\/{\rm )}$\times\mathbb{R}^n$ 
if and only if the following two hold:   
$$ 
Q\left(LJ_\Phi\right)\cong Q(x),\; 
\frac{\partial LJ_\Phi}{\partial x}(0,0)\ne 0.       
$$ 
\end{enumerate}
\end{thm}
Here, $Q(f_1,\ldots, f_\ell) $ stands for Mather's local algebra for function-germs 
$f_1, \ldots, f_\ell$.      For Mather's local algebra, see \cite{mather4, wall}.       
Theorems \ref{theorem 3} (resp., Theorem \ref{theorem 4}) 
is used as a criterion of Whitney umbrella 
(resp., (swallowtail)$\times\mathbb{R}^{n-1}$).   
Theorems \ref{theorem 3} and \ref{theorem 4} 
are connected by the following simple lemma.   
\begin{lmm}\label{lemma 1}
\quad 
For a normalized Legendrian map-germ 
$\Phi: (\mathbb{R}\times\mathbb{R}^n,(0,0))\to (\mathbb{R}^2\times\mathbb{R}^n,(0,0))$, 
$$
LJ_{\Phi}(x,y)=(-1)^{n+1}
\frac{\frac{\partial \Phi_2}{\partial x}(x,y)}{\nu_1(x,y)}.
$$
Here 
$\nu_{\Phi}(x,y)  =  \nu_1(x,y)\frac{\partial}{\partial X_1}+\nu_2(x,y)
\frac{\partial}{\partial X_2}+\cdots +\nu_{n+2}(x,y)\frac{\partial}{\partial X_{n+2}}$.  
\end{lmm}
In the case $n=1$, Lemma \ref{lemma 1} with its proof can be found in \cite{nishimurapacific}.   
The proof given in \cite{nishimurapacific} works well in general case.       

\section{Applications of calculus correspondence} \label{section 3}
In order to show that the calculus correspondence is significant and useful, we give two applications 
of Theorem \ref{theorem 1}.   
\begin{prp}\label{application}
Let $\Phi: (\mathbb{R}\times\mathbb{R}^n,(0,0))\to (\mathbb{R}^2\times\mathbb{R}^n,(0,0))$ 
be given by 
$$
\Phi(x,y)=\left( ax^4+x^2\sum_{i=1}^n b_iy_i, cx^3+x\sum_{i=1}^n d_iy_i, y\right)
\quad (a,b_i,c,d_i\in \mathbb{R}) 
$$    
where $y=(y_1, \ldots, y_n)$.     
Then, the following two are equivalent.   
\begin{enumerate}
\item The given 
$\Phi$ is a {\rm (}swallowtail\/{\rm )}$\times\mathbb{R}^{n-1}$, 
that is, it is $\mathcal{A}$-equivalent to the normal form of 
{\rm (}swallowtail\/{\rm )}$\times\mathbb{R}^{n-1}$ 
which is the following 
$$
(x,y)\mapsto \left(3x^4+x^2y_1, -4x^3-2xy_1, y\right).
$$ 
\item The equality $2ad_i=3b_ic$ holds for any $i$ $(1\le i\le n)$ 
and there exists an $i$ $(1\le i\le n)$ such that $ab_icd_i\ne 0$ is satisfied. 
\end{enumerate}    
\end{prp}     
\textit{Proof of Proposition \ref{application}.}\qquad 
Suppose that $\Phi$ is a (swallowtail)$\times\mathbb{R}^{n-1}$.      
Then, since $\Phi$ is Legendrian, there exists a unit normal vector field  
$$
\nu_{\Phi}(x,y)=(\nu_1(x,y), \nu_{2}(x,y), \ldots , \nu_{n+2}(x,y))
$$ 
such that the following two hold: 
\begin{equation}\label{x}
\nu_1(x,y)(4ax^3+2x\sum_{i=1}^n b_iy_i)+\nu_2(x,y)(3cx^2+\sum_{i=1}^n d_iy_i)  =  0, 
\end{equation}
\begin{equation}\label{y}
\nu_1(x,y)b_ix^2+\nu_2(x,y)d_ix+\nu_{2+i}(x,y)=0\quad (1\le i\le n).
\end{equation}
Since $\Phi$ is a (swallowtail)$\times\mathbb{R}^{n-1}$, 
$\nu_2(0,0)$ (resp., $\nu_{2+i}(0,0)$) must be zero by the equality (\ref{x}) (resp., (\ref{y})).     
Since $\nu_{\Phi}(0,0)$ is a unit vector, $\nu_1(0,0)$ must be $\pm 1$.    
It is clear that the given $\Phi$ satisfies the first and the second conditions of Definition \ref{normalized}.            
Thus, $\Phi$ is a normalized (swallowtail)$\times\mathbb{R}^{n-1}$.     
By Theorem \ref{theorem 1}, 
$$
\mathcal{D}(\Phi)(x,y)=\left(4ax^3+2x\sum_{i=1}^n b_iy_i, 3cx^2+\sum_{i=1}^n d_iy_i, y \right)
$$ 
is a (Whitney umbrella)$\times\mathbb{R}^{n-1}$ of pedal unfolding type.     
Since $\mathcal{D}(\Phi)$ is of pedal unfolding type, we have that $2ad_i=3b_ic$ for any $i$ $(1\le i\le n)$.        
Since $\mathcal{D}(\Phi)$ is a (Whitney umbrella)$\times\mathbb{R}^{n-1}$, 
there must exist an $i$ $(1\le i\le n)$ such that 
$b_ic\ne 0$.    
Two conditions $2ad_i=3b_ic$ and $b_ic\ne 0$ imply that $ab_icd_i\ne 0$.    
 \par 
 Conversely, suppose that 
the equality $2ad_i=3b_ic$ holds for any $i$ $(1\le i\le n)$ 
and there exists an $i$ $(1\le i\le n)$ such that $ab_icd_i\ne 0$ is satisfied. 
%
Then, $\mathcal{D}(\Phi)$ is a (Whitney umbrella)$\times\mathbb{R}^{n-1}$ of pedal unfolding.         
Thus, by Theorem \ref{theorem 1}, $\Phi=\mathcal{I}(\mathcal{D}(\Phi))$ is a normalized 
(swallowtail)$\times\mathbb{R}^{n-1}$.   
\hfill\qquad $\Box$ 
\par 
\medskip 
As another application of Theorem \ref{theorem 1}, we give an alternative proof of 
Arnol'd's observation given in \cite{arnolddevelopable} (for Arnol'd's observation, see also \cite{ishikawa2}),        
Namely, we show the following:   
\begin{observation}\label{observation}
Let $\gamma : (\mathbb{R}, 0)\to (\mathbb{R}^3,0)$ be the space curve 
given by $\gamma(x)=(x^4,x^3,x^2)$.         Then, the tangent developable of 
$\gamma$, which is the following, is a swallowtail.   
$$
\Phi(x,y)=\left(x^4, x^3, x^2\right)+y\left(4x^2, 3x, 2\right).   
$$  
\end{observation}
\textit{Proof of Observation \ref{observation}.}\qquad 
Put $\widetilde{y}=x^2+2y$.    
Then, $\Phi$ is $\mathcal{R}$-equivalent to 
$\widetilde{\Phi}(x,\widetilde{y})
=\left(-x^4+2x^2\widetilde{y}, -\frac{1}{2}x^3+\frac{3}{2}x\widetilde{y}, \widetilde{y}
\right)$.     
It is easily seen that 
$\mathcal{D}(\widetilde{\Phi})$ 
is a Whitney umbrella of pedal unfolding type.      
Thus, by Theorem \ref{theorem 1}, $\widetilde{\Phi}=\mathcal{I}(\mathcal{D}(\widetilde{\Phi}))$ 
is a normalized swallowtail.   
\hfill\qquad $\Box$
%
%

\section{Questions around calculus correspondences} \label{section 5}

The following question is a multi-parameter version of the question posed in \cite{nishimurapacific}.   
\begin{question}\label{question 2}
\begin{enumerate}
\item Let $\varphi_1, \varphi_2: (\mathbb{R}\times\mathbb{R}^n,(0,0))\to 
(\mathbb{R}^2\times\mathbb{R}^n,(0,0))$ 
be two map-germs of pedal unfolding type.   
Suppose that $\varphi_1$ is $\mathcal{A}$-equivalent to $\varphi_2$.   
Is $\mathcal{I}(\varphi_1)$ necessarily $\mathcal{A}$-equivalent to $\mathcal{I}(\varphi_2)$ ?   
\item Let $\Phi_1, \Phi_2: (\mathbb{R}\times\mathbb{R}^n,(0,0)) \to 
(\mathbb{R}^2\times\mathbb{R}^n,(0,0))$ be two normalized Legendrian  map-germs.   
Suppose that $\Phi_1$ is $\mathcal{A}$-equivalent to $\Phi_2$.   
Is $\mathcal{D}(\Phi_1)$ necessarily $\mathcal{A}$-equivalent to $\mathcal{D}(\Phi_2)$ ?   
\end{enumerate}
\end{question}
Question \ref{question 2} seems to be difficult to solve completely in general.            
In the following two subsections,  we discuss special cases of Question \ref{question 2}.

\subsection{$S_k$ type singularities and Legendrian $S_k$ type singularities} 

Recall that a map-germ $f: (\mathbb{R}^2,0)\to (\mathbb{R}^3,0)$ is of $S_k$ type 
if $f$ is $\mathcal{A}$-equivalent to the map-germ 
$f_{k,\pm}(x,y)=\left(x\left(x^2\pm y^{k+1}\right), x^2\pm y^{k+1}, y\right)$ 
(Definition \ref{definition S_k}).      Since the map-germ $f_{k,\pm}$ is of 
pedal unfolding type, the following map-germ (which is $\mathcal{I}(f_{k, \pm})$) is normalized Legendrian map-germ by 
Proposition \ref{proposition 4}.          
$$
F_{k,\pm}(x,y)=\left(\frac{1}{4}x^4\pm \frac{1}{2}x^2y^{k+1}, \frac{1}{3}x^3\pm xy^{k+1}, y\right).   
$$
The Legendrian map-germ $\mathcal{I}(f_{k, \pm})$ is called   
the {\it normal form of Legendrian $S_k$ type} and any Legendrian map-germ $\mathcal{A}$-equivalent to 
$\mathcal{I}(f_{k, \pm})$ is said to be {\it of Legendrian $S_k$ type}.  
\begin{question}\label{question 3}
\begin{enumerate}
\item Let $\varphi: (\mathbb{R}\times\mathbb{R},(0,0))\to 
(\mathbb{R}^2\times\mathbb{R},(0,0))$ 
be a map-germ of pedal unfolding type.   
Suppose that $\varphi$ is of $S_k$ type.     
Is $\mathcal{I}(\varphi)$ necessarily of Legendrian $S_k$ type ?   
\item Let $\Phi:(\mathbb{R}\times\mathbb{R},(0,0)) \to 
(\mathbb{R}^2\times\mathbb{R},(0,0))$ be a normalized Legendrian  map-germ.   
Suppose that $\Phi$ is of Legendrian $S_k$ type.          
Is $\mathcal{D}(\Phi)$ necessarily of $S_k$ type ?   
\end{enumerate}
\end{question}
In the case $k=0$, both $f_{0,+}, f_{0,-}$ are $\mathcal{A}$-equivalent to 
the normal form of Whitney umbrella, and both $F_{0,+}, F_{0,-}$  are 
$\mathcal{A}$-equivalent to the normal form of swallowtail (namely, the map-germ 
$(x,y)\mapsto (3x^4+x^2y, -4x^3-2xy, y)$).              
In this case, we have the calculus correspondence by Theorem \ref{theorem 1}.       
\par 
In the case $k=1$, $f_{1,+}$ (resp., $F_{1,+}$) is 
not $\mathcal{A}$-equivalent to $f_{1, -}$ (resp., $F_{1,-}$).        
It is known that 
only the map-germs of $S_1$ type are $\mathcal{A}_e$-codimension one singularities of mono-germs 
from the plane to the 3-space (for $\mathcal{A}_e$-codimension, see \cite{wall} and 
for the classification of $\mathcal{A}_e$-codimension one singularities 
$(\mathbb{R}^2,0)\to (\mathbb{R}^3,0)$, 
see \cite{chenmatumoto, coopermondatique, mond}).            
Theorem \ref{theorem 3} can be applied as a criterion of $S_1$ singularities.       
On the other hand, criteria of Legendrian $S_1$ singularities have been obtained by Izumiya-Saji-Takahashi 
(\cite{izumiyasajitakahashi}).         
Thus, by replacing Saji-Umehara-Yamada criterion (Theorem \ref{theorem 4}) 
with Izumiya-Saji-Takahashi criteria given in \cite{izumiyasajitakahashi}, the proof of Theorem \ref{theorem 1} is 
expected to work well to show calculus correspondence between $S_1$ singularities of pedal unfolding type and 
normalized Legendrian $S_1$ singularities.    
\par 
Next, we discuss the case $k\ge 2$.       
Even in this case, Theorem \ref{theorem 3} can be applied as a criterion of $S_k$ singularities.       
However, there seems to be no criteria for Legendrian $S_k$ singularities in the case $k\ge 2$.           
Hence, it seems that we cannot expect an analogy of the proof of Theorem \ref{theorem 1}.    
        
\subsection{Legendrian $A_k$ type singularities}
\begin{dfn}[\cite{sajiumeharayamada}]
{\rm 
Let $k, n$ be non-negative integers such that $k\le n+1$.     
\begin{enumerate}
\item    The map-germ 
$G_k: (\mathbb{R}\times\mathbb{R}^{n},(0,0))\to (\mathbb{R}^2\times\mathbb{R}^{n},(0,0))$ given by 
$$
G_k(x,y)=
\left((k+1)x^{k+2}+\sum_{j=1}^{k-1} jx^{j+1}y_j, -(k+2)x^{k+1}-\sum_{j=1}^{k-1} (j+1)x^j y_j, y\right)
$$
is called the {\it normal form of Legendrian $A_{k+1}$ type}, where 
$(x,y)=(x, y_1, \ldots, y_n)$.    
\item    A map-germ $\Phi: (\mathbb{R}\times\mathbb{R}^n,(0,0))\to 
(\mathbb{R}^2\times\mathbb{R}^n,(0,0))$ is said to be 
{\it of} {\it Legendrian $A_{k+1}$ type\/} if $\Phi$ is $\mathcal{A}$-equivalent to $G_k$.    
\end{enumerate}
}
\end{dfn}
Note that the image of $G_k$ is the envelope of the following one parameter family of 
hyperplanes.       By this reason, $G_k$ is called the normal form of $A_{k+1}$ type. 
$$
\left\{\left(X_1, X_2, Y_1, \ldots, Y_n\right)\; |\; 
x^{k+2}+Y_{k-1}x^k+\cdots +Y_1x^2+X_2x+X_1=0.\right\}    
$$
For the normal form of Legendrian $A_{k+1}$ type, we have 
$$
\mathcal{D}(G_k)(x,y)=\left(n(x,y)p(x,y), p(x,y), y\right), 
$$
where $n(x,y)=-x$ and 
$p(x,y)=-(k+2)(k+1)x^k-\sum_{j=1}^{k-1}j(j+1)x^{j-1}y_j$.       
Since, $p(0,0)=0$ and $\frac{\partial n}{\partial x}(0,0)\ne 0$, $\mathcal{D}(G_k)$ is of pedal unfolding type.    
Therefore, $G_k=\mathcal{I}(\mathcal{D}(G_k))$ is normalized Legendrian by 
Proposition \ref{proposition 4}.          
\begin{question}\label{question Gaffney}
\begin{enumerate}
\item Let $\varphi: (\mathbb{R}\times\mathbb{R}^n,(0,0))\to 
(\mathbb{R}^2\times\mathbb{R}^n,(0,0))$ 
be a map-germ of pedal unfolding type.   
Suppose that $\varphi$ is $\mathcal{A}$-equivalent to $\mathcal{D}(G_k)$.      
Is $\mathcal{I}(\varphi)$ necessarily of Legendrian $A_{k+1}$ type ?   
\item Let $\Phi:(\mathbb{R}\times\mathbb{R}^n,(0,0)) \to 
(\mathbb{R}^2\times\mathbb{R}^n,(0,0))$ be a normalized Legendrian  map-germ.   
Suppose that $\Phi$ is of Legendrian $A_{k+1}$ type.          
Is $\mathcal{D}(\Phi)$ necessarily $\mathcal{A}$-equivalent to $\mathcal{D}(G_k)$ ?   
\end{enumerate}
\end{question}
\noindent 
Question \ref{question Gaffney} was asked by G.~Ishikawa (\cite{ishikawa}), and independently 
by T.~Gaffney during AMS Spring Western Section Meeting at the University of Hawaii (2012).     
It is easily seen that $G_1$ is non-singular, $G_2$ is the normal form of 
(cusp)$\times\mathbb{R}^n$ and $G_3$ is the normal form of 
(swallowtail)$\times\mathbb{R}^{n-1}$.           Thus, in the case $k=0, 1, 2$, 
Proposition \ref{proposition 5}, Theorem \ref{theorem 2} and Theorem \ref{theorem 1} 
are the affirmative answers to Question \ref{question Gaffney} respectively.   
\par 
Therefore, Question \ref{question Gaffney} asks essentially the case $k\ge 3$.         
Even in this case, there is a criterion of Legendrian $A_{k+1}$ singularities (Theorem \ref{A_k criterion}).        
However, there seems to be 
no criteria for the $\mathcal{A}$-equivalence class of $\mathcal{D}(G_k)$ $(k\ge 3)$.   
Hence, it seems that we cannot expect an analogy of the proof of Theorem \ref{theorem 1}.    
\begin{thm}[\cite{sajiumeharayamada}] \label{A_k criterion}
For a normalized Legendrian map-germ 
$\Phi: (\mathbb{R}\times\mathbb{R}^n,(0,0))$ 
$\to (\mathbb{R}^2\times\mathbb{R}^n,(0,0))$, 
$\Phi$ is of Legendrian $A_{k+1}$ type if and only if the following two hold:    
$$ 
Q\left(LJ_\Phi, \frac{\partial LJ_\Phi}{\partial x}, \ldots, \frac{\partial^{k-1} LJ_\Phi}{\partial x^{k-1}}\right)
\cong Q(x,y_1, \ldots, y_{k-1}), \; 
\frac{\partial^k LJ_\Phi}{\partial x^k}(0,0)\ne 0.    
$$
\end{thm}
\subsection*{Acknowledgement}
The author 
would like to express his sincere gratitude to the referee who was kind enough 
to give a simple and clear proof of Proposition \ref{proposition 1} and make valuable suggestions.   

\end{document}